\documentclass[11pt,british]{amsart}

\newif\ifpdf
\ifx\pdfoutput\undefined\pdffalse\else\pdfoutput=1\pdftrue\fi\newif\pdf
\ifpdf\relax\else
\usepackage[T1]{fontenc}
\fi

\usepackage{babel,eucal,url,amssymb,stmaryrd,
enumerate,amscd,
}

\usepackage[pagebackref]{hyperref}

\def\CC {{\mathbb C}}     
\def\RR {{\mathbb R}}

\theoremstyle{plain}
\newtheorem{thm}{Theorem}[section]
\newtheorem{lem}[thm]{Lemma}
\newtheorem{pro}[thm]{Proposition}
\newtheorem{co}[thm]{Corollary}

\theoremstyle{definition}
\newtheorem{defn}[thm]{Definition}

\theoremstyle{remark}
\newtheorem{rem}[thm]{Remark}

\newtheorem*{ack}{Acknowledgements}

\newcommand{\Gtwo}{\ifmmode{{\rm G}_2}\else{${\rm G}_2$}\fi}

\def\mk {\mathfrak}

\date{\today}
\begin{document}

\title[Hyper-ParaHermitian  manifolds with torsion]%
{Hyper-ParaHermitian  manifolds with torsion}

\author{Stefan Ivanov}
\address[Ivanov]{University of Sofia "St. Kl. Ohridski"\\
Faculty of Mathematics and Informatics,\\ Blvd. James Bourchier
5,\\ 1164 Sofia, Bulgaria,}
\email{ivanovsp@fmi.uni-sofia.bg}

\author{Vasil Tsanov}
\address[Tsanov]{University of Sofia "St. Kl. Ohridski"\\
Faculty of Mathematics and Informatics,\\ Blvd. James Bourchier
5,\\ 1164 Sofia, Bulgaria,}
\email{tsanov@fmi.uni-sofia.bg}

\author{Simeon Zamkovoy}
\address[Zamkovoy]{University of Sofia "St. Kl. Ohridski"\\
Faculty of Mathematics and Informatics,\\
Blvd. James Bourchier 5,\\
1164 Sofia, Bulgaria}
\email{}

\begin{abstract}
Necessary and sufficient conditions for the existence of a
hyper-parahermitian connection with totally skew-symmetric torsion
(HPKT-structure) are presented. It is shown that any
HPKT-structure is locally generated by a real (potential)
function. An invariant first order differential operator is
defined on any almost hyper-paracomplex manifold showing that it
is two-step nilpotent exactly when the almost hyper-paracomplex
structure is integrable. A local HPKT-potential is expressed in
terms of this operator. Examples of (locally) invariant
HPKT-structures with closed as well as non-closed torsion 3-form
on a class of (locally) homogeneous hyperparacomplex manifolds
(some of them compact) are constructed.

\end{abstract}

\maketitle
\setcounter{tocdepth}{2}
\tableofcontents

\section{Introduction}
We study the geometry of structures on a differentiable manifold
related to the algebra of paraquaternions together with a
naturally associated metric which is necessarily of neutral
signature. This structure leads to the notion of (almost)
hyper-paracomplex and hyper-parahermitian manifolds in dimensions
divisible by four. These structures are also attractive in
theoretical physic since some of them play a role in string theory
\cite{OV,hul,Bar,Hull} and integrable systems \cite{D1}.

Hyper-parahermitian geometry may be interpreted as the indefinite
analog of hyper-hermitian geometry, but there are important
differences. We provide hyper-parahermitian versions of
many local and some global results for hyper-hermitian manifolds,
specially we adopt the hypercomplex constructions of
\cite{HP,GP,BS}  (but see also \cite{S3,MCS,MS,Ver}).

We treat integrable almost hyper-parahermitian structures,
which admit compatible linear connections with totally-skew
symmetric torsion, briefly HPKT-structure. It is known that in
dimension 4, the conformal structure of neutral signature
determined by a hyper-paracomplex structure is necessarily
anti-self-dual \cite{AG,hul,IZ}. We show that the corresponding
conformal hyper-parahermitian structure is an HPKT-structure.
In higher dimensions, we find necessary and sufficient conditions for
the existence of a HPKT-structure in terms of the exterior derivative
of the three K\"ahler forms. We give a holomorphic characterization
and show uniqueness of the HPKT connection.

To illustrate the subtleties of HPKT we use some homogeneous
examples and their compact factors found in \cite{ITZ}. In
particular, we show the existence of an invariant HPKT-structure
with closed torsion 3-form on the simple Lie groups $SU(m,m-1),
m>1$, associated to the biinvariant Killing-Cartan form of
neutral signature on $SU(m,m-1)$.
In contrast, the HPKT-structures for the hyper-paracomplex
structures on the simple Lie groups $SL(2m-1,\mathbb R),m>1$
obtained in \cite{ITZ} have no compatible biinvariant metric.
They may be associated to (a class of) invariant
metrics of neutral signature, which however have non-closed
torsion forms.

We show that any HPKT-structure is locally generated by a real
(potential) function following the ideas developed in \cite{BS}.
To this end, using Salamon's idea from the
quaternionic case (see \cite{S3}), we define an invariant first
order differential operator $D$, the hyper-paracomplex operator,
on an almost hyper-paracomplex manifold and we show that it is
2-step nilpotent exactly when the almost hyper-paracomplex
structure is integrable. Then we obtain the local existence of
HPKT-potential by proving the local $D$-exactness of
certain $D$-closed 2-forms.

\begin{ack}
The research is partially supported by Contract MM 809/1998 with
the Ministry of Science and Education of Bulgaria, Contracts
586/2002 and 35/2003 with the University of Sofia "St. Kl.
Ohridski". S.I. is a  member of the EDGE, Research Training
Network HPRN-CT-2000-00101, supported by the European Human
Potential Programme. We would like to thank to T. Gramtchev for
his valuable comments.
\end{ack}

\section{HyperparaK\"ahler connection with torsion}

Both quaternions $H$ and paraquaternions $\tilde H$ are real
Clifford algebras, $H=C(2,0),\quad \tilde H=C(1,1)\cong C(0,2)$.
In other words, the algebra $\tilde H$ of paraquaternions is
generated by the unity $1$ and the generators $J_1,J_2,J_3$
satisfying the \emph{paraquaternionic identities},
\begin{equation}\label{par1}
J_1^2=J_2^2=-J_3^2=1,\qquad J_1J_2=-J_2J_1=J_3.
\end{equation}

We recall the notion of almost hyper-paracomplex manifold
introduced by Libermann \cite{Lib}. An \emph{almost quaternionic
structure of the second kind}  on a smooth manifold consists of
two almost product structures $J_1,J_2$ and an almost complex
structure $J_3$ which mutually anti-commute, i.e. these structures
satisfy the paraquaternionic identities \eqref{par1}. Such a
structure is also called \emph{complex product structure}
\cite{ASal,An}.

An \emph{almost hyper-paracomplex structure} on a 4n-dimensional
manifold $M$ is a triple $\tilde H=(J_a), a=1,2,3$, where
$J_{\alpha}$,$\alpha = 1,2$ are almost paracomplex structures
$J_{\alpha}:TM\rightarrow TM$, and $J_{3}:TM\rightarrow TM$ is an
almost complex structure, satisfying the paraquaternionic
identities \eqref{par1}. We note that on an almost
hyper-paracomplex manifold there is actually a 2-sheeted
hyperboloid worth of almost complex structures:
$S^2_1(-1)=\{c_1J_1+c_2J_2+c_3J_3 : c_1^2+c_2^2-c_3^2=-1\}$ and a
1-sheeted hyperboloid worth of almost paracomplex structures:
$S^2_1(1)=\{b_1J_1+b_2J_2+b_3J_3 : b_1^2+b_2^2-b_3^2=1\}$.

When each $J_a,a=1,2,3$ is an integrable structure, $\tilde H$ is
said to be a \emph{hyper-paracomplex structure} on $M$. Such a
structure is also called sometimes \emph{pseudo-hyper-complex}
\cite{D1}. Any hyper-paracomplex structure admits a  unique
torsion-free connection $\nabla^{CP}$ preserving $J_1,J_2,J_3$
\cite{ASal,An} called \emph{the complex product connection}.

The Nijenhuis tensor $N_a$ of $J_a$ is defined by
\begin{equation}\label{nu}
N_{\alpha}(X,Y)=[J_{\alpha}X,J_{\alpha}Y]+J_{\alpha}^{2}[X,Y]-
J_{\alpha}[J_{\alpha}X,Y]
- J_{\alpha}[X,J_{\alpha}Y].
\end{equation}

It is well known that the structure $J_a$ is integrable if and
only if the corresponding Nijenhuis tensor $N_a$ vanishes,
$N_a=0$.

In fact an almost hyper-paracomplex structure is hyper-paracomplex
if and only if any two of the three structures $J_a,a=1,2,3$ are
integrable due to the existence of a linear identity between the
three  Nijenhuis tensors \cite{IZ,BVuk0}. In this case all almost
complex structures of the two-sheeted hyperboloid  $S^2_1(-1)$ as
well as all paracomplex structures of the one-sheeted hyperboloid
$S^2_1(1)$ are integrable.

A \emph{hyperparahermitian metric} is a pseudo Riemannian metric
which is compatible with the (almost) hyperparacomplex structure
$\tilde H=(J_a), a=1,2,3$ in the sense that the metric $g$ is
skew-symmetric with respect  to each $J_a, a=1,2,3$, i.e.
\begin{equation}\label{par2}
g(J_1.,J_1.)=g(J_2.,J_2.)=-g(J_3.,J_3.)=-g(.,.).
\end{equation}
The metric $g$ is necessarily of neutral signature (2n,2n). Such a structure
is
called \emph{(almost) hyper-paraHermitian structure}.

Let $F_a$ be the K\"ahler form associated with the structure
$(g,J_a), a=1,2,3$:$$F_a=g(.,J_a.).$$ The corresponding Lee form
is defined by $\theta_a=-\delta F_a\circ J_a^3$. In particular,
$$\theta_a(X)=\sum_{i=1}^{2n}dF_a(e_i,J_ae_i,J_a^2X),$$ for an
orthonormal $J_a$-adapted basis $\{e_1,\dots,
e_{2n},J_ae_1,\dots,J_ae_{2n}\}$.

If on a hyper-paraHermitian manifold there exists an admissible
basis $(\tilde H)$ such that each $J_a, a=1,2,3$ is parallel with
respect to the Levi-Civita connection or equivalently the three
K\"ahler forms are cosed, $dF_a=0$ then the manifold is called
\emph{hyper-paraK\"ahler}. Such manifolds are also called
\emph{hypersymplectic} \cite{Hit}, \emph{neutral hyper-K\"ahler}
\cite{Kam,FPed}. The equivalent characterization is that the
holonomy group of $g$ is contained in $Sp(n,\mathbb R)$ if $ n\ge
2$ \cite{Vuk}.

For $n=1$ an (local) almost hyper-paracomplex structure is the
same as oriented neutral conformal structure
\cite{D1,GRio,Vuk,BVuk0}. The existence of a (local)
hyper-paracomplex structure is a strong condition since the
integrability of the (local) almost hyper-paracomplex structure
implies that the corresponding neutral conformal structure is
anti-self-dual \cite{AG,hul,IZ}. The necessary and sufficient condition
for the integrability of an (local) almost hyper-paracomplex
structure in dimension four is the coincidence of the three Lee
forms, $\theta_1=\theta_2=\theta_3$ \cite{IZ}.

We  use the following notations: For any r-form $\omega$ we define
$J_a\omega(X_1,\dots,X_r):=(-1)^r\omega(J_aX_1,\dots,J_aX_r),
a=1,2,3$ and the operators
$d_{\alpha}\omega:=-J_{\alpha}dJ_{\alpha}\omega, \alpha=1,2,\\
d_3\omega:=(-1)^rJ_3dJ_3\omega$. In particular
$d_aF_a=-dF_a(J_a.,J_a.,J_a), \quad a=1,2,3$.

We consider the (para) complex operators
\begin{gather*}
\partial_{\alpha}=\frac{1}{2}\left(d+\epsilon d_{\alpha}\right),
\quad \bar\partial_{\alpha}=\frac{1}{2}\left(d-\epsilon
d_{\alpha}\right),\quad \epsilon^2=1, \quad \alpha=1,2 \\
\partial_3=\frac{1}{2}\left(d+i d_3\right),
\quad \bar\partial_3=\frac{1}{2}\left(d-i d_3\right),\quad i^2=-1.
\end{gather*}
In particular, a complex function $f=u+iv$ is holomorphic with
respect to the complex structure $J_3$ iff $\bar\partial_3 f=0$
while a paracomplex function $h=u+\epsilon v$ is paraholomorphic
with respect to the paracomplex structure $J_{\alpha}, \alpha=1,2$
iff $\bar\partial_{\alpha} h=0$.
\begin{defn}
A hyperparahermitian metric $g$ is \emph{hyperparaK\"ahler
with torsion}
(briefly HPKT) if there exists a linear connection $\nabla$
preserving the hyperparacomplex structure whose torsion tensor
$T^{\nabla}$ is totally skew-symmetric i.e.
\begin{equation}\label{fff}
\nabla g= \nabla J_1=\nabla J_2=\nabla J_3=0, \quad
T^{\nabla}(X,Y,Z):=g(T^{\nabla}(X,Y),Z)=-T^{\nabla}(X,Z,Y).
\end{equation}
If the torsion 3-form $T^{\nabla}$ is closed, $dT^{\nabla}=0$,
then the HPKT-metric is called \emph{strong HPKT metric}.
\end{defn}
A connection satisfying condition (\ref{fff}) will be called
briefly \emph{HPKT-connection}.
\begin{thm}
Let $(M,g,J_1,J_2,J_3)$ be a  hyperparahermitian manifold. The
following conditions are equivalent:
\begin{enumerate}
\item{$(M,g,J_1,J_2,J_3)$
admits a HPKT-connection;}
\item{The following equalities hold
\begin{equation}\label{hp1}
d_1F_1=d_2F_2=d_3F_3.
\end{equation}
In this case the HPKT connection is uniquely determined by the
torsion
\begin{equation}\label{hp2}
T^{\nabla}=d_1F_1=d_2F_2=d_3F_3.
\end{equation}}
\end{enumerate}
In particular, the three Lee forms coincide,
$\theta_1=\theta_2=\theta_3=tr_gT^{\nabla}$.
\end{thm}
\begin{proof}
The required connection is the unique Bismut connection determined
by Gauduchon \cite{Ga1} (see also \cite{FI}) in the hermitian case
and by Ivanov-Zamkovoy \cite{IZ} in the parahermitian case due to
the compatibility condition \eqref{hp1}. Take the trace in
\eqref{hp2} to get the last identity.
\end{proof}
It is known that in dimension four any hyperparahermitian metric
is anti-self-dual \cite{AG,hul,IZ}. The proof of
Theorem~6.2 in \cite{IZ} leads to
\begin{pro}
Any hyperparahermitian metric on a hyperparacomplex 4-manifold is
HPKT. In particular, the Ricci two-forms of the HPKT-connection
all vanish.
\end{pro}
\section{Homogeneous examples}
A non-trivial class of examples for the
differential geometric entities defined in the previous
section is provided by certain
left-invariant HPKT-structures on (semi)
simple Lie groups which were found in \cite{ITZ}.
For convenience we
reproduce here the explicit description of $J_2$ and $J_3$. We
define the (para-) complex structures on Lie
algebras and then interpret them as homogeneous
almost (para-) complex
structures on the corresponding simply connected
Lie groups. For brevity we shall
sometimes abuse notation and proper definitions, by indicating
only the Lie algebras.

\subsection{HPKT-structure on $SU(m,m-1)$}
The most important example from \cite{ITZ}is the group $SU(m,m-1)$,
where the biinvariant Killing form is the neutral HPKT metric.

On the simple Lie algebra $su(m,m-1)$ (of dimension $4m(m-1)$)
we define a scalar product
\begin{gather}\label{scal}
B(X,Y) \doteq \frac{1}{2}tr(XY),\quad X,Y \in su(m,m-1)
\end{gather}
Obviously $B$ is proportional to the Killing form and
defines a biinvariant, neutral pseudoriemannian metric on
$SU(m,m-1)$.
Next we produce a convenient $B$ - orthonormal base of the
Lie algebra $su(m,m-1)$. As usual, we
denote by $E_j^k\in gl(n)$ the matrix with entry 1 at
the intersection of the j-th row and the k-th column and 0
elsewhere. We fix the range of indices
\begin{gather}\label{inde}
j = 1,\dots,m -1 ,\qquad j < k < 2m-j.
\end{gather}

Let $\mk{z}$ be the subspace (abelian subalgebra)
of $su(m,m-1)$ generated by the elements
\begin{gather*}
i(E^j_j + E^{2m-j}_{2m-j} - 2E^m_m),\quad j = 1,\dots,m-1.
\end{gather*}
Let $Z^1,\dots.Z^{m-1}$ be any orthonormal base of $\mk{z}$,
with respect to the scalar product\footnote{ $B$ is
obviously negative definite on $\mk{z}$.} $B$. We define
\begin{gather}\label{sell5}
X^j \doteq i(E^j_j - E^{2m-j}_{2m-j});
\quad
Y^j \doteq E_j^{2m-j} + E_{2m-j}^j:
\quad
W^j\doteq i(E_j^{2m-j} - E_{2m-j}^j).
\end{gather}
\begin{equation}\label{sell15}
\begin{split}
U = U_j^k \doteq &
\begin{cases}
E_j^k - E_k^j& \text{ if }\quad  j< k \leq m; \\ E_k^{2m - j} -
E_{2m-j}^k & \text{ if }\quad m < k < 2m-j.
\end{cases}\\
V = V_j^k \doteq &
\begin{cases}
i(E_j^k + E_k^j)& \text{ if }\quad  j< k \leq m;\\ i(E_k^{2m - j}
+ E_{2m-j}^k) & \text{ if }\quad m < k < 2m-j.
\end{cases}\\
S = S_j^k \doteq &
\begin{cases}
E_k^{2m - j}+ E_{2m-j}^k& \text{ if }\quad  j < k \leq m; \\ =
E_j^k - E_k^j & \text{ if }\quad m < k < 2m-j.
\end{cases}\\
T = T_j^k \doteq &
\begin{cases}
i(E_k^{2m - j}- E_{2m-j}^k)& \text{ if }\quad  j < k \leq m; \\
 i(E_k^j - E_j^k)& \text{ if }\quad m < k < 2m-j.
\end{cases}
\end{split}
\end{equation}
The invariant vector fields (generated by)
$X^j, Y^j, W^j, Z^j, U_j^k, V_j^k, S_j^k, T_j^k$ give
a base of the tangent bundle of $SU(m,m-1)$. We define an almost
hyperparacomplex structure by
\begin{align}\label{sell3}
J_3(Z^j) \doteq X^j;\quad &J_3(Y^j) \doteq W^j;\quad &J_2(Z^j)
\doteq W^j,\quad &J_2(X^j) \doteq Y^j;\\ \notag J_3(U_j^k) \doteq
V_j^k;\quad &J_3(S_j^k) \doteq T_j^k;\quad &J_2(U_j^k) \doteq
T_j^k; \quad &J_2(V_j^k) \doteq S_j^k
\end{align}
It is shown in \cite{ITZ} that the structure \eqref{sell3} is an
integrable hyper-paracomplex structure on $SU(m,m-1), m>1$ which
is compatible with the biinvariant (Killing-Cartan) form of
neutral signature $B$.

Now we observe, that the above construction gives also a strong
HPKT-structure on $SU(m,m-1)$. The HPKT-connection is the
left-invariant connection $\nabla$, defined by postulating all
left-invariant vector fields to be parallel.

The torsion of the above connection is the Lie bracket
and the torsion tensor $T^{\nabla}(X,Y,Z)=-B([X,Y],Z)$ is a
closed 3-form (due to the Jacobi identity). So, we have
a strong HPKT-structure on $SU(m,m-1)$ which is flat. The
compatible neutral Killing-Cartan metric is Einstein.

Simple Lie groups admit cocompact lattices \cite{Bor},
say $\Gamma$. Hence, we obtain a HPKT interpretation of the
result proved in \cite{ITZ}.

\begin{thm}\cite{ITZ}
The compact manifolds $SU(m,m-1)/\Gamma$ admit invariant, flat,
strong HPKT-structures. The neutral HPKT-metric is a non-flat
Einstein metric induced by the Killing-Cartan form.
\end{thm}

\begin{rem}
The above procedure can be applied to the group
$(SL(2m-1,\CC))^{\RR}$ (see \cite{ITZ}). Thus we obtain invariant
strong and flat HPKT-structures on the compact manifolds
$(SL(2m-1,\CC))^{\RR}/\Gamma$.
\end{rem}

\subsubsection{A non-strong HPKT-structure on $SU(2,1)$}
We equipped the 8-dimensional simple Lie group $SU(2,1)$ with a
strong and flat left-invariant HPKT-structure induced by the
left-invariant hyper-paracomplex structure \eqref{sell3} and the
Killing-Cartan form. We show below that a similar\footnote{
We choose the simplest one in a notational sense. Obviously, any
neutral metric on the Lie algebra $\mk{g}$, which is compatible
with a paraquaternionic structure, gives a left invariant metric
on the
corresponding Lie group $G$. However there is only one
biinvariant metric on a simple group.}
hyper-paracomplex structure supports left-invariant HPKT-structure
which is not strong and not flat . To be precise, we consider the
following base on $su(2,1)$
\begin{align}\label{selln5}
Z \doteq &i(E^1_1 + E^3_3 - 2 E^2_2);\quad &X \doteq i(E^1_1 -
E^3_3);\\ \notag
W \doteq &i(E_1^3 - E_3^1);\quad &Y \doteq E_1^3
+ E_3^1; \\ \notag
U  \doteq & E_1^2 - E_2^1;\quad &V \doteq
i(E_1^2 + E_2^1);\\ \notag
S \doteq & E_2^3+ E_3^2; \quad &T
\doteq i(E_2^3- E_3^2).
\end{align}
A hyper-paracomplex structure on $SU(2,1)$ is given by
\begin{align}\label{sela}
J_3(Z) \doteq X;\quad &J_3(Y) \doteq W;\quad &J_2(Z) \doteq
W,\quad &J_2(X) \doteq Y;\\ \notag J_3(U) \doteq V;\quad &J_3(S)
\doteq T;\quad &J_2(U) \doteq T; \quad &J_2(V) \doteq S,
\end{align}
We claim that the neutral  metric $g$ determined by the following
orthonormal basis
\begin{gather}\label{met2}
g(Z,Z)=g(X,X)=g(U,U)=g(V,V)=1,\\ \nonumber
g(Y,Y)=g(W,W)=g(S,S)=g(T,T)=-1
\end{gather}
is a non strong left-invariant HPKT-metric on $SU(2,1)$ with
respect to the left-invariant hyper-paracomplex structure
\eqref{sela}. We denote the 1-form dual to a vector field via the
neutral metric \eqref{met2} by the same letter. We calculate
\begin{gather*}
T^{\nabla}=d_1F_1=d_2F_2=d_3F_3\\\nonumber =2X\wedge Y\wedge
W-X\wedge U\wedge V+X\wedge S\wedge T+Y\wedge U\wedge
S\\\nonumber-Y\wedge V\wedge T +W\wedge U\wedge T+W\wedge V\wedge
S-Z\wedge U\wedge V-Z\wedge S\wedge T;\\\nonumber dT^{\nabla}=-4
U\wedge V\wedge S\wedge T\not=0.
\end{gather*}
Our claim is proved.

\subsection{HPKT-structure on $SL(2m-1,\mathbb R)$}

Consider the simple Lie group $SL(2m-1,\mathbb R)$ with the almost
hyper-paracomplex structure \eqref{sell3} applied to the following
base\footnote{The range of indices is as in
formula (\ref{inde}).} of $sl(2m-1,\mathbb R)$:
\begin{align}\label{sell14}
Z^j \doteq &E^j_j + E^{2m-j}_{2m-j} - 2 E^m_m; \quad &W^j \doteq
E^j_j - E^{2m-j}_{2m-j};\\ \notag
X^j \doteq &E_j^{2m-j} -
E_{2m-j}^j;\quad &Y^j \doteq E_j^{2m-j} + E_{2m-j}^j;\\ \notag
U_j^k \doteq  &E_j^k - E_k^j;\quad &V_j^k \doteq E_k^{2m - j} -
E^k_{2m-j};\\ \notag S_j^k \doteq &E^{2m - j}_k + E^k_{2m-j};
\quad, &T_j^k \doteq E_j^k + E_k^j.
\end{align}
The structure \eqref{sell3} is an integrable left-invariant
hyper-paracomplex structure on $SL(2m-1,\RR), m>1$ \cite{ITZ}.

A left-invariant neutral metric $g$, determined by the following
orthonormal basis of $sl(2m-1,\RR)$
\begin{gather}\label{met1}
g(Z^j,Z^j)=g(X^j,X^j)=g(U^k_j,U^k_j)=g(V^k_j,V^k_j)=1,\\ \nonumber
g(Y^j,Y^j)=g(W^j,W^j)=g(S^k_j,S^k_j)=g(T^k_j,T^k_j)=-1
\end{gather}
is an HPKT-metric on $sl(2m-1,\RR),m>1$ with respect to the
hyper-paracomplex structure \eqref{sell3}, which has not closed
torsion 3-form. Indeed, denoting the 1-form dual to a vector field
via the neutral metric \eqref{met1} by the same letter, we
calculate
\begin{gather*}
T^{\nabla}=d_1F_1=d_2F_2=d_3F_3\\\nonumber =2X^j\wedge Y^j\wedge
W^j-X^j\wedge U^k_j\wedge V^k_j+X^j\wedge S^k_j\wedge
T^k_j+Y^j\wedge U^k_j\wedge S^k_j\\\nonumber-Y^j\wedge V^k_j\wedge
T^k_j +W^j\wedge U^k_j\wedge T^k_j+W^j\wedge V^k_j\wedge
S^k_j+Z^j\wedge U^k_j\wedge T^k_j-Z^j\wedge V^k_j\wedge
S^k_j;\\\nonumber dT^{\nabla}=-8 U^k_j\wedge V^k_j\wedge
S^k_j\wedge T^k_j\not=0.
\end{gather*}

 The groups $SL(2m-1,\mathbb R)$ admit cocompact lattices
 \cite{Bor}, say $\Gamma$. Thus, we arrive at a
 HPKT-extension of the results in \cite{ITZ}.

\begin{thm}
The compact manifolds $SL(2m-1,\mathbb R)/\Gamma$ admit invariant
HPKT-structure which are not strong.
\end{thm}

\subsection{HPKT-structure on $2\mathbb R\oplus sl(2,\mathbb C)$}
We consider the following base on $2\RR \oplus sl(2,\CC)$:
\begin{gather}\label{ccs21}
Z \doteq
\begin{bmatrix}
1 & 0\\ 0 & 1
\end{bmatrix}
X \doteq
\begin{bmatrix}
0 & 1\\ -1 & 0
\end{bmatrix} ;
Y \doteq
\begin{bmatrix}
 0 & 1\\
1 & 0
\end{bmatrix};
W \doteq
\begin{bmatrix}
1 & 0\\ 0 &-1
\end{bmatrix}\\\nonumber
U \doteq iZ;\quad V \doteq iX;\quad S \doteq iY;\quad T \doteq iW.
\end{gather}
We define an almost hyper-paracomplex structure on the Lie algebra
$2\RR \oplus sl(2,\CC) \cong 2\RR \oplus so(3,1)$ by \eqref{sela}
using the base \eqref{ccs21} . It is easy to check that this
structure  is integrable.

We claim that the neutral  metric $g$ determined by the
orthonormal basis \eqref{met2} is a strong left-invariant
HPKT-metric on the simply connected Lie group $G$ associated to
the Lie algebras $2\RR \oplus sl(2,\CC) \cong 2\RR \oplus
so(3,1)$ with respect to the left-invariant
hyper-paracomplex structure \eqref{sela}.
To prove the claim, we denote the 1-form dual to a
vector field via the neutral metric \eqref{met2} by the same
letter. We obtain
\begin{equation*}
T^{\nabla}=d_1F_1=d_2F_2=d_3F_3=S\wedge dS-Y\wedge dY; \qquad
dT^{\nabla}=0
\end{equation*}
which proves our claim.

Let $\Delta
\subset SL(2,\CC)$ be a cocompact discrete subgroup and let
$\Gamma=\mathbb Z\times \mathbb Z\times \Delta\subset G$. We obtain
\begin{thm}
The compact manifold $G/\Gamma$ admits an invariant strong
HPKT-structure.
\end{thm}

\section{Characterization of HPKT-structures}
In this section, we characterize HPKT-structures in terms of the
existence of holomorphic objects. We use ideas from the definite
(HKT) case described in \cite{HP,GP} and find other
compact examples.
\subsection{Holomorphic characterization}
We recall that the space of paracomplex (1,0)-vectors (resp.
(0,1)-vectors) with respect to the paracomplex structure
$J_{\alpha},\alpha=1,2$ is spanned by paracomplex vectors of type
$X+\epsilon J_{\alpha}X$ (resp. $X-\epsilon J_{\alpha}X$) and the
space of complex (1,0)-vectors (resp. (0,1)-vectors) of the
complex structure $J_3$ is spanned as usual by complex vectors of
type $X-iJ_3X$ (resp. $X+iJ_3X$).

It is easy to check that
\begin{description}
\item{}the 2-form $F_2-\epsilon F_3$ is of type (2,0)
while the 2-form $F_2+\epsilon F_3$ is of type (0,2) with respect
to the paracomplex structure $J_1$;
\item{}the 2-form $F_3+\epsilon F_1$ is of type (2,0)
while the 2-form $F_3-\epsilon F_1$ is of type (0,2) with respect
to the paracomplex structure $J_2$;
\item{}the 2-form $F_1-i F_2$ is of type (2,0)
while the 2-form $F_1+ i F_2$ is of type (0,2) with respect to the
complex structure $J_3$;
\end{description}
\begin{pro}\label{hpro1}
Let $(M,g,J_a, a=1,2,3)$ be a hyperparahermitian manifold. The
following conditions are equivalent:
\begin{description}
\item[a] $(M,g,J_a, a=1,2,3)$ is a PHKT manifolld;
\item[b] $\partial_1(F_2-\epsilon F_3)=0$ or equivalently
$\bar\partial_1(F_2+\epsilon F_3)=0$;
\item[c]$\partial_2(F_3+\epsilon F_1)=0$ or equivalently
$\bar\partial_2(F_3-\epsilon F_1)=0$;
\item[d] $\partial_3(F_1-i F_2)=0$ or equivalently
$\bar\partial_3(F_1+ i F_2)=0$;
\end{description}
\end{pro}
\begin{proof}
We have $$\partial_1 (F_2 - \epsilon F_3)=\bar{\partial_1} (F_2 +
\epsilon
F_3)=\frac{1}{2}(dF_2-d_1F_3)-\frac{\epsilon}{2}(dF_3-d_1F_2).$$
Therefore $\partial_1 (F_2 - \epsilon F_3)=0$, when the real and
imaginary parts both vanishes. We calculate
\begin{gather*}d_1F_3=-J_1dJ_1F_3=-J_1d(F_3 \circ J_1)=-J_1dF_3=(dF_3 \circ
J_1)=(dF_3 \circ J_3J_2)\\ =-J_3(dF_3 \circ
J_2)=J_2J_3dF_3=J_2J_3dJ_3F_3=J_2d_3F_3. \end{gather*} On the
other hand $$dF_2=-d(F_2 \circ J_2)=-J_{2}^2dJ_2F_2=J_2d_2F_2.$$
Consequently, the condition $d_1F_3=dF_2$ is equivalent to the
condition $d_2F_2=d_3F_3$. Therefore the Bismut connection of the
parahermitian structure $(g,J_2)$ coincides with Bismut connection
of the hermitian structure $(g,J_3)$. Since $J_1=J_3J_2$ then
$J_1$ is parallel with respect to the common connection $\nabla$.
Therefore $\nabla$ is the Bismut connection for $(g,J_1)$ which
proves the equivalence of a) and b). In a similar way one
completes the proof.
\end{proof}
Proposition~\ref{hpro1} implies that the HPKT condition is not
preserved by a generic conformal transformation of the metric
provided the dimension is at least eight.

In the proof of Proposition~\ref{hpro1}, we also derive
\begin{co}
Suppose $F_1,F_2$ and $F_3$ are the K\"ahler forms of a
hyperparahermitian structure. Then the hyperparahermitian
structure is HPKT-structure if and only if
$$d_aF_b=\delta_{ab}T^{\nabla}-\epsilon_{abc}F_c,$$ where
$\delta_{ab}$ is the Kroneker symbol and $\epsilon_{abc}$ is the
totally skew-symmetric Levi-Civita symbol.
\end{co}
\begin{thm}\label{thex}
Let $(M,J_a, a=1,2,3)$ be a hyperparacomplex manifold. Then any
one of the following three conditions implies the forth:
\begin{enumerate}
\item{} $F_2+ \epsilon F_3$ is a $(0,2)$-form with respect to
$J_1$ such that $\bar{\partial_1} (F_2 + \epsilon F_3)=0$ and
$F_2(X,J_2Y)=g(X,Y)$ is a symmetric non-degenerate bilinear form
of neutral signature;
\item{} $F_3- \epsilon F_1$ is a $(0,2)$-form with respect to $J_2$
such that $\bar{\partial_2} (F_3 - \epsilon F_1)=0$ and
$F_1(X,J_1Y)=g(X,Y)$ is a symmetric non-degenerate bilinear form
of neutral signature;
\item{} $F_1+ iF_2$ is a $(0,2)$-form with respect to $J_3$ such
that $\bar{\partial_3} (F_1 + iF_2)=0$ and $F_3(X,J_3Y)=g(X,Y)$ is
a symmetric non-degenerate bilinear form of neutral signature;
\item{}The structure $(g,J_a, a=1,2,3)$ is a $PHKT$ structure.
\end{enumerate}
\end{thm}
\begin{proof}
In view of the Proposition~\ref{hpro1}, it suffices to prove that
the metric $g$ is hyperparahermitian.

Using the fact that $F_2+ \epsilon F_3$ is of type $(0,2)$ with
respect to $J_1$. Since $X+ \epsilon J_1X$ is of type $(1,0)$ with
respect to $J_1$, $(F_2+ \epsilon F_3)(X+ \epsilon J_1X,Y)=0$, for
any vectors $X,Y$. It is equivalent to the identity
$F_3(X,Y)=-F_2(J_1X,Y)$. Then
$$F(X,J_3Y)=-F_2(J_1X,J_3Y)=-F_2(J_1X,J_1J_2Y)=-F_2(X,J_2Y)=-g(X,Y).$$
So $F_3(J_3X,J_3Y)=F_3(X,Y)$ and $g$ is hermitian with respect to
$J_3$. Since the metric is parahermitian with respect to $J_2$ and
$J_1=J_3J_2$, $g$ is parahermitian with respect to $J_1$.

Similarly, one get the other assertions.
\end{proof}

\subsection{HPKT structures on compact Nilmanifolds}\label{sec1}

In this section we construct further examples of
homogeneous HPKT-structures, now on some (compact) nilmanifolds.

Let $\{X_1,\dots,X_{2n},Y_1,\dots,Y_{2n},Z\}$ be a basis for
$\mathbb R^{4n+1}$. Define commutators by: $[X_j,Y_j]=Z$, all
others being zero. These commutators give $\mathbb R^{4n+1}$
the structure of the \emph{Heisenberg Lie agebra} $\mk{h}_{2n}$.
Let
$\mathbb R^3$ be the three dimensional abelian algebra. The direct
sum $\mk{n}= \mk{h}_{2n}\oplus \mathbb R^3$ is a
2-step nilpotent algebra
whose center is four-dimensional. Fix a basis $\{E_1,E_2,E_3\}$
for $\mathbb R^3$ and consider the following endomorphisms of
$\mk{n}$:
\begin{gather}\label{hais}\nonumber
J_2 : X_{2j-1} \rightarrow Y_{2j}, X_{2j} \rightarrow Y_{2j-1}
\quad Z \rightarrow E_2,\quad E_1\rightarrow -E_3; \\  J_3 :
X_{2j-1} \rightarrow X_{2j}, \quad Y_{2j-1} \rightarrow Y_{2j}
\quad Z \rightarrow E_1,\quad E_2\rightarrow E_3; \\ \nonumber
J_2^2=-J_3^2= identity, \quad J_1=J_3J_2.
\end{gather}
Clearly $J_2J_3=-J_3J_2$. The almost complex structure $J_3$
satisfies the identity $[J_3.,J_3.]=[.,.]$ which implies that it
is an Abelian almost complex structure on $\mk{n}$
in the sense of
\cite{BMM} and in particular integrable.
The next notion seems to be useful
\begin{defn}
The almost paracomplex
structure $J_2$ is said to be \emph{Abelian} if the following
identity $[J_2.,J_2.]=-[.,.]$ holds.
\end{defn}
Applying \eqref{nu} it is easy to check that any Abelian almost
paracomplex structure has vanishing Nijenuis tensor and therefore
is integrable. It is easy to verify that the almost paracomplex
structure $J_2$ is Abelian on $\mk{n}$. Consequently, the almost
paracomplex structure $J_1$ is also Abelian.  Hence, the structure
$J_a,a=1,2,3$ is a left invariant hyperparacomplex structure on
the simply connected Lie group $N$ whose Lie algebra is $\mk{n}$.
Consider the invariant metric $g$ on $N$ for which the basis
$\{X_j,Y_j,Z,E_a\}$ is orthonormal and
$g(X_j,X_j)=g(Z,Z)=g(E_1,E_1)=1, \quad
g(Y_j,Y_j)=g(E_2,E_2)=g(E_3,E_3)=-1$. Clearly, the structure
$(g,J_a,a=1,2,3)$ is a left invariant hyperparahermitian structure
on $N$ which turns out to be a HPKT since any left invariant
(2,0)-form with respect to the complex structure $J_3$ is
$\partial_3$-closed due to a result of Salamon \cite{saln} and
Proposition~\ref{hpro1}. Because $N$ is isomorphic to the
product $H_{2n}\times \mathbb R^3$ of the Heisenbrg group $H_{2n}$
and the Abelian group $\mathbb R^3$ we have:
\begin{co}
Let $\Gamma$ be a cocompact lattice in the Heisenberg group
$H_{2n}$ and $\mathbb Z^3$ a lattice in $\mathbb R^3$. The compact
Nilmanifold $N/(\Gamma\times\mathbb Z^3)$ admits an invariant
HPKT-structure.
\end{co}

\subsection{HPKT-structure on $(H_{2n}\times \widetilde{SL(2,\mathbb
R)})/\Gamma$} Based on the above computations, we can also see
that there is a left-invariant HPKT-structure on the product of
$4n+1$-dimensional Heisenberg group $H_{2n}$ and the universal
cover $\widetilde{SL(2,\mathbb R})$ of the simple Lie group
$SL(2,\mathbb R)$. The Lie algebra $sl(2,\mathbb R)$ has a basis
$\{E_1,E_2,E_3\}$ with non-zero brackets given by
\\ \centerline{$[E_1,E_2]=E_3, \quad [E_2,E_3]=-E_1, \quad
[E_3,E_1]=E_2$.} The construction of the HPKT-structure on the
product $H_{2n}\times\widetilde{SL(2,\mathbb R)}$ is the same as
those in Section~\ref{sec1} taking $E_1,E_2,E_3$ to be the
generators of $sl(2,\mathbb R)$. The integrability of the
(non-abelian) almost hyper-paracomplex structure defined by
\eqref{hais} as well as the HPKT-compatibility conditions
\eqref{hp1} can be checked directly using the commutators of the
left-invariant vector fields. Denote the left-invariant 1-forms
dual to the left invariant vector fields via the metric  by the
same letters to get
\begin{equation*}
T^{\nabla}=d_1F_1=d_2F_2=d_3F_3=dZ\wedge Z; \qquad
dT^{\nabla}=dZ\wedge dZ\not=0.
\end{equation*}
The last equalities imply that the HPKT-structure is not strong.

Let $\Gamma_1$ be a cocompact lattice in the Heisenberg group
$H_{2n}$. The universal cover $\widetilde{SL(2,\mathbb R)}$ of the
Lie group ${SL(2,\mathbb R)}$ admits a discrete subgroup
$\Gamma_2$ such that the quotient space $\widetilde{(SL(2,\mathbb
R)}/\Gamma_2)$ is a compact 3-manifold  \cite{Mil,RV,Scot}. Such a
space has to be Seifert fibre space \cite{Scot} and all the
quotients are classified in \cite{RV}. We obtain
\begin{co}
The compact manifold $(H_{2n}\times \widetilde{SL(2,\mathbb
R)})/(\Gamma_1\times\Gamma_2)$ admits an invariant non-strong
HPKT-structure.
\end{co}

\section{Potential theory}
It is well known that a K\"ahler metric is locally generated by a
potential ie a real function $\mu$ satisfying $F_3=-dd_3\mu$.
Similarly a paraK\"aler metric is locally generated by a
potential, ie a real function $\nu$ satisfying $F_1=dd_1\nu$
\cite{Rash}.

A function $\mu$ is a potential function for a hyperparaK\"ahler
manifold $(M,g,J_a)$ if the K\"ahler forms are equal to
\begin{equation}\label{kf1}
F_a=J_a^2dd_a\mu, \quad d_a\mu=-J_a^3d\mu.
\end{equation}

In this section, we seek a function that generates all K\"ahler
forms of a HPKT-manifold.

The definition of the operators $d_a$, paraquaternionic identities
\eqref{par1}, the compatibility conditions \eqref{par2} and
\eqref{kf1} imply
\begin{gather}\label{poten1}
d_1d_2\mu=-d_1J_2d\mu=J_1dJ_1J_2d\mu=J_1dd_3\mu=-J_1F_3=-d_2d_1\mu=dd_3\mu;\\
\nonumber
d_2d_3\mu=d_2J_3d\mu=J_2dJ_1d\mu=-J_2dd_1\mu=-J_2F_1=-d_3d_2\mu=-dd_1\mu;\\
\nonumber
d_3d_1\mu=-d_3J_1d\mu=J_3dJ_3J_1d\mu=J_3dd_2\mu=J_3F_2=-d_1d_3\mu=-dd_2\mu;
\nonumber
\end{gather}
We generalize this concept to HPKT-manifold.
\begin{defn}
Let $(M,g,J_a)$ be a HPKT-structure with K\"ahler forms $F_1,F_2$
and $F_3$. A possibly locally defined function $\mu$ is a
potential function for the HPKT structure if
\begin{gather}\label{poten2}
F_1=\frac{1}{2}\left(dd_1-d_2d_3\right)\mu, \quad
F_2=\frac{1}{2}\left(dd_2-d_3d_1\right)\mu, \quad
F_3=-\frac{1}{2}\left(dd_3+d_1d_2\right)\mu.
\end{gather}
\end{defn}
In fact any one of the above identities implies the others due to
the next
\begin{thm}\label{po1}
Let $(M,g,J_a)$ be a HPKT-structure with K\"ahler forms $F_1,F_2$
and $F_3$. Let $\nabla^{CP}$ be the complex product connection. A
possibly locally defined function $\mu$ is a potential function
for the HPKT structure if any one of the following identities hold
\begin{gather}
F_1=\frac{1}{2}\left(dd_1-d_2d_3\right)\mu, \label{po2}\\
F_2=\frac{1}{2}\left(dd_2-d_3d_1\right)\mu, \label{po3}\\
F_3=-\frac{1}{2}\left(dd_3+d_1d_2\right)\mu, \label{po4}\\
g=\frac{1}{2}\left(1-J_1-J_2+J_3\right)(\nabla^{CP})^2\mu.\label{po5}
\end{gather}
The torsion 3-form $T^{\nabla}$ is given by $T^{\nabla}= -
\frac{1}{2}d_1d_2d_3\mu$.
\end{thm}
\begin{proof}
We calculate, using the fact that the complex product connection
is a torsion-free and  preserves the hyper-paracomplex structure,
that
\begin{gather*}
dd_3\mu(X,Y)=-(\nabla^{CP}_Xd\mu)J_3Y+(\nabla^{CP}_Yd\mu)J_3X;\\
d_1d_2\mu(X,Y)=J_1dd_3\mu(X,Y)=(\nabla^{CP}_{J_1X}d\mu)J_2Y-
(\nabla^{CP}_{J_1Y}d\mu)J_2X;\\
g(X,Y)=-F_3(X,J_3Y)=\frac{1}{2}\left(dd_3+d_1d_2\right)\mu=
\frac{1}{2}\left(1-J_1-J_2+J_3\right)(\nabla^{CP})^2\mu
\end{gather*}
Thus, the equivalence of \eqref{po4} and \eqref{po5} is proved.

Similarly one can get the equivalence of \eqref{po2} and
\eqref{po5} as well as the equivalence between \eqref{po3} and
\eqref{po5}.

The formula for the torsion is a consequence of \eqref{poten2} and
\eqref{hp2}.
\end{proof}
\begin{rem}
In the context of a potential, an HPKT structure is
hyper-paraK\"ahler if and only if the potential function $\mu$
satisfies any of the following four identities
\begin{gather*}
dd_1\mu=-d_2d_3\mu,\\ dd_2\mu=-d_3d_1\mu,\\ dd_3\mu=d_1d_2\mu,\\
(1+J_1+J_2+J_3)(\nabla^{CP})^2\mu=0.
\end{gather*}
\end{rem}
\begin{co}
Let $(M,g,J_a)$ be a HPKT-structure with K\"ahler forms $F_1,F_2$
and $F_3$. A possibly locally defined function $\mu$ is a
potential function for the HPKT structure if any one of the
following identities hold
\begin{gather}
F_2-\epsilon F_3= - 2
\partial_1J_2\bar\partial_1\mu; \label{poten3} \\
F_3+\epsilon F_1= - 2\partial_2J_3\bar\partial_2\mu;
\label{poten4}\\ F_1-iF_2= - 2\partial_3J_1\bar\partial_3\mu.
\label{poten5}
\end{gather}
\end{co}
\begin{proof}
Due to \eqref{poten1} and the definition of the operators
$\partial_a, \quad \bar\partial_a$, we have
\begin{equation*}
F_1-iF_2=\frac{1}{2}\left(dd_1-d_2d_3 -idd_2+id_3d_1\right)= -
2\partial_3J_1\bar\partial_3\mu
\end{equation*}
The other assertions follow in a similar way.
\end{proof}
\subsection{HPKT potential in dimension 4}
Now we give a hyperbolic version  of the existence of the HKT potentials of Banos
and Swann \cite{BS}. We apply Theorem~\ref{po1} to the 4-dimensional
case to prove the existence of a local HPKT potential for any
HPKT metric.
\begin{co}
Let $g$ be an HPKT metric on a 4-dimensional hyper-paracomplex
manifold and let $\theta$ be  1-form defined by the complex
product connection via $\nabla^{CP}g=\theta\otimes g$.

A function $\mu$ is an HPKT potential for $g$ if and only if it is
solution of the hyperbolic equation
\begin{equation*}
\triangle\mu -d\mu(\theta^{\sharp}) +2=0,
\end{equation*}
where $\triangle$ is the hyperbolic Laplacian of the neutral
metric $g$.

In particular, any HPKT metric on a 4-dimensional
hyper-paracomplex manifold admits locally a potential.
\end{co}
\begin{proof}
Let $A=\nabla^{CP}-\nabla^g$, where $\nabla^g$ denote the
Levi-Civita connection of $g$. The tensor $A$ is symmetric,
$A(X,Y)=A(Y,X)$ since both connections are torsion-free. We also
have $$\theta(X)g(Y,Z)=-g(A(X,Y),Z)-g(A(X,Z),Y).$$ Solving for
$A$, we obtain
$$g(A(X,Y),Z)=\frac{1}{2}\left(\theta(Z)g(X,Y)-\theta(X)(g(Y,Z)-\theta(Y)g(Y,Z)\right).
$$ In particular, if $X$ is a (local) unit vector field, then
$$g(A(X,X),Y)=\frac{1}{2}\theta(Y)-\theta(X)g(X,Y)$$ and
$$g(A(X,X)-A(J_1X,J_1X)-A(J_2X,J_2X)+A(J_3X,J_3X),Y)=\theta(Y)$$
for all $Y$.

The metric $g$ is the unique hyper-parahermitian metric satisfying
$g(X,X)=1$. Therefore $\mu$ is a HPKT potential if and only if
$$\frac{1}{2}\left(1-J_1-J_2+J_3\right)(\nabla^{CP})^2(X,X)=1,$$
that is $$Trace_g(\nabla^{CP}d\mu)=2.$$ Note that the hyperbolic
Laplacian $\triangle\mu$ is by definition
$-Trace_g(\nabla^gd\mu)$. Thus $\mu$ is a HPKT potential for $g$
if and only if
$$-\triangle\mu+d\mu(A(X,X)-A(J_1X,J_1X)-A(J_2X,J_2X)+A(J_3X,J_3X))=2.
$$ The local existence of HPKT potentials now follows from the
general theory for the (ultra) hyperbolic Laplace operator (see
e.g. \cite{H1} and the references therein).
\end{proof}

\subsection{HPKT-potential in dimension $4n\ge 8$} Here we
demonstrate the existence of a local potential for any HPKT-metric
(the HKT case was done by Banos and Swann \cite{BS}).

The crucial step is the construction on any almost hyperparacomplex
manifold of an $GL(n,\tilde H)$-invariant first order differential
operator D, which is the hyperbolic version of the hypercomplex
(quaternionic) operator of Salamon (see \cite{S3}). The operator $D$ is
two-step nilpotent if and only if the structure is hyperparacomplex.
Then we obtain the existence of a local HPKT-potential in terms of the
operator D.

The element $$\dag=-J_1\otimes J_1-J_2\otimes J_2+J_3\otimes J_3$$
is independent of the choice of the basis $\{J_1,J_2,J_3\}$ and
acts naturally on $\Lambda^2$ with $\dag^2=2\dag +3$. The
eigenspace decomposition
\begin{equation}\label{op1}
\Lambda^2=\{\dag=3\}\oplus\{\dag=-1\}
\end{equation}
is a paraquaternionic invariant in the sense that it is preserved
by $GL(n,\tilde H)Sp(1,\mathbb R)$ and therefore it is a
hyper-paracomplex invariant preserving by $GL(n,\tilde H)$.

\subsection{The hyper-paracomplex differential.} Studying the
action of $GL(n,\tilde H)$ on the bundle $\Lambda^k$ we consider
the subbundle
\begin{equation*}
A^k=\sum_{I\in
S^2_1(-1)}\left(\Lambda^{k,0}_I\oplus\Lambda^{0,k}_I\right).
\end{equation*}
It is not difficult to see that
\begin{equation*}
A^k=\sum_{P\in
S^2_1(1)}\left(\Lambda^{k,0}_P\oplus\Lambda^{0,k}_P\right).
\end{equation*}
Indeed,  any 2-form $\omega\in\Lambda^2$ decomposes according to \eqref{op1}
\begin{gather}\label{op2}
\omega(X,Y)=\frac{1}{4}\left\{3\omega(X,Y)+\omega(J_1X,J_1Y)+
\omega(J_2X,J_2Y)-\omega(J_3X,J_3Y)\right\}\\ \nonumber
+\frac{1}{4}\left\{\omega(X,Y)-\omega(J_1X,J_1Y)-
\omega(J_2X,J_2Y)+\omega(J_3X,J_3Y)\right\}\\ \nonumber
\end{gather}
For example
\begin{equation*}
A^2=\Lambda^{2,0}_{J_3}\oplus\Lambda^{0,2}_{J_3}\oplus
A^{1,1}_{J_3}=\Lambda^{2,0}_{J_2}\oplus\Lambda^{0,2}_{J_2}\oplus
A^{1,1}_{J_2},
\end{equation*}
where
\begin{gather*}
A^{1,1}_{J_3}=\{\omega\in\Lambda^2 : J_3\omega=\omega \quad {\rm
and} \quad J_2\omega=\omega\},\\
A^{1,1}_{J_2}=\{\omega\in\Lambda^2 : J_2\omega=-\omega \quad {\rm
and}\quad J_3\omega=-\omega\}.
\end{gather*}
If $g$ is an hyper-parahermitian metric then the K\"ahler form
$F_3$, (resp. $F_2$) is a smooth section of $A^{1,1}_{J_3}$,
(resp. $A^{1,1}_{J_2}$) and conversely any smooth section $F_3$ of
$A^{1,1}_{J_3}$, (resp. $F_2$ of $A^{1,1}_{J_2}$) defines an
(possibly degenerate) hyper-parahermitian metric $g=-F_3(.,J_3.)$,
(resp. $g=F_2(.,.J_2)$. We will call such a form  \emph{a
hyper-paracomplex (1,1)-form}.

There is a projection $\eta : \Lambda^k\rightarrow A^k$ whose
kernel is the subbundle
\begin{gather*}
B^k=\bigcap_{I\in
S^2_1(-1)}\left(\Lambda^{k-1,1}_I\oplus\Lambda^{k-2,2}_I\oplus
\dots\oplus \Lambda^{1,k-1}_I \right)=\\
=\bigcap_{P\in
S^2_1(1)}\left(\Lambda^{k-1,1}_P\oplus\Lambda^{k-2,2}_P\oplus
\dots\oplus \Lambda^{1,k-1}_P \right).
\end{gather*}
In particular, the two eigenspaces of the operator $\dag$ are
related with $A^2,B^2$ as follows $$A^2=\{\dag=-1\}, \quad
B^2=\{\dag=3\}. $$ The projections $\omega^{A^2}$ and
$\omega^{B^2}$ are given by
\begin{gather*}
\omega^{A^2}(X,Y)=\frac{1}{4}\left\{3\omega(X,Y)+\omega(J_1X,J_1Y)+
\omega(J_2X,J_2Y)-\omega(J_3X,J_3Y)\right\},\\
\omega^{B^2}(X,Y)=\frac{1}{4}\left\{\omega(X,Y)-\omega(J_1X,J_1Y)-
\omega(J_2X,J_2Y)+\omega(J_3X,J_3Y)\right\}.
\end{gather*}

We define the \emph{hyper-paracomplex differential}
\begin{equation*}
D : A^k\rightarrow A^{k+1}
\end{equation*}
simply by composition of the projection $\eta$ end the exterior
differential $d$:
\begin{equation*}
D=\eta\circ d.
\end{equation*}
For example, if $\omega$ is a 1-form then
\begin{gather}\label{dhp1}
D\omega=(d\omega)^{2,0}_{J_3}+(d\omega)^{0,2}_{J_3}+\frac{1}{2}\left((d\omega)^{1,1}_{J_3}
+J_2(d\omega)^{1,1}_{J_3}\right)=\\ \nonumber =
(d\omega)^{2,0}_{J_2}+(d\omega)^{0,2}_{J_2}+\frac{1}{2}\left((d\omega)^{1,1}_{J_2}
-J_3(d\omega)^{1,1}_{J_2}\right.)
\end{gather}
\begin{thm}\label{int}
An almost hyper-paracomplex structure is integrable if and only if
$D^2=0$.
\end{thm}
\begin{proof} The condition $D^2=0$
is equivalent to the assertion that the exterior derivative of a
2-form of type (1,1) relative to all $I\in S^2_1(-1)$ and all
$P\in S^2_1(1)$ has no (0,3)+(3,0)-component relative to any $I\in
S^2_1(-1)$ and to any $P\in S^2_1(1)$. The latter condition holds
on a hyper-paracomplex manifold since all almost complex
structures of the two-sheeted hyperboloid $S^2_1(-1)$ as well as
all paracomplex structures of the one-sheeted hyperboloid
$S^2_1(1)$ are integrable due to the existence of a linear
identity between their  Nijenhuis tensors \cite{IZ,BVuk0}.

To prove the converse, let $\Omega$ be (1,1)-form with respect to
the almost complex structure $J_3, \Omega\in\Lambda^{1,1}_{J_3}$.
The 2-form $C$ defined by $C(X,Y)=\Omega(X,Y)-\Omega(J_1X,J_1Y)$
belongs to $B^2$. The condition $DC=0$ is equivalent to the
relation $$dC(X,Y,Z)=dC(IX,IY,Z)+dC(IX,Y,IZ)+dC(X,IY,IZ)$$ for all
$I \in S^2_1(-1)$ and
$$-dC(X,Y,Z)=dC(PX,PY,Z)+dC(PX,Y,PZ)+dC(X,PY,PZ)$$ for all $P \in
S^2_1(1)$. In particular
\begin{gather}\label{91}
-dC(J_3X,J_3Y,J_3Z)=\\ \nonumber
dC(J_1X,J_1Y,J_3Z)+dC(J_1X,J_3Y,J_1Z)+dC(J_3X,J_1Y,J_1Z).
\end{gather}
Let $\nabla^{CP}$ be a complex product connection,
$\nabla^{CP}J_a=0$, with torsion $T^{CP}$. Then the  Nijenhuis
tensors are related with $T^{CP}$ as follows
\begin{equation}\label{new1}
N_a=-T^{CP}(J_a.,J_a.)-J_a^2T^{CP}(.,.)+J_aT^{CP}(J_a.,.)+J_aT^{CP}(.,J_a.).
\end{equation}
Use \eqref{new1} to express the exterior derivative of a 2-form as
\begin{gather}\label{92}
dC(X,Y,Z)= \nabla^{CP}C(X;Y,Z)+ \nabla^{CP}C(Y;Z,X)+
\nabla^{CP}C(Z;X,Y)\\ \nonumber
+C(T^{CP}(X,Y),Z)+C(T^{CP}(Y,Z),X)+C(T^{CP}(Z,X),Y).
\end{gather}
Insert \eqref{92} into \eqref{91} and use \eqref{new1} to get
\begin{gather}\label{new2}
\Omega(J_3N_2(J_3X,J_3Y),Z)
+\Omega(J_3N_2(J_3Y,J_3Z),X)+\Omega(J_3N_2(J_3Z,J_3X),Y)+\\
\nonumber \Omega(J_1N_2(J_3X,J_3Y),J_2Z)
+\Omega(J_1N_2(J_3Y,J_3Z),J_2X)+\Omega(J_1N_2(J_3Z,J_3X),J_2Y)=0
\end{gather}
valid for any (1,1)-form with respect to $J_3$. In particular,
take $\Omega=Z\wedge J_3Z$, we get from \eqref{new2} that $N_2=0$.

Similarly, we obtain $N_1=0$. Hence, the almost hyper-paracompex
structure is integrable \cite{IZ}.
\end{proof}
Further we need two lemmas
\begin{lem}\label{int2}
Let $F\in A^{1,1}_{J_3}$ be a non-degenerate hyperparacomplex
(1,1)-form on a hyper-paracomplex manifold $(M,J_1,J_2,J_3)$. The
metric $g=-F(.,J_3.)$ is HPKT metric if and only if $F$ is
$D$-closed, $DF=0$.

Such a form is called a HPKT-form
\end{lem}
\begin{proof}
Suppose that $g$ is HPKT. For any complex structure $I \in
S^2_1(-1)$ (paracomplex strucrure $P \in S^2_1(1)$) the form
$dF_I$ ($dF_P$) has type $(2,1)+(1,2)$ with respect to the complex
structure $I$ (paracomplex structure $P$). But since
$d_1F_1=d_2F_2=d_3F_3$ we deduce that $dF_I$ ($dF_P$) has type
$(2,1)+(1,2)$ with respect to the two paracomplex structures and
the complex structure: $dF_I \in B^3$ ($dF_P \in B^3$) that is
$DF_I=0$ ($DF_P=0$). Since $F_I=F$ ($F_P=F$) we obtain the result.
\\
Suppose now that $DF=0$. The condition $F \in A^{(1,1)}_{J_3}$
also reads as
\begin{equation}\label{93}
F(X,Y)=F(J_3X,J_3Y)=F(J_1X,J_1Y), F(J_3X,J_1Y)=F(J_1X,J_3Y).
\end{equation}
Use the torsion free complex product connection, \eqref{92},
\eqref{93} and \eqref{91}, to get
\begin{gather*}
-dF(J_3X,J_3Y,J_3Z)=\nabla^{CP}F(J_3X;Y,Z)+
\nabla^{CP}F(J_3Y;Z,X)+ \nabla^{CP}F(J_3Z;X,Y)\\
+2\nabla^{CP}F(J_1X;J_1Y,J_3Z)+ 2\nabla^{CP}F(J_1Y;J_1Z,J_3X)+
2\nabla^{CP}F(J_1Z;J_1X,J_3Y).
\end{gather*}
Consequently
\begin{gather*}
-dF(J_3X,J_3Y,J_3Z)=\\\nabla^{CP}F(J_1X;J_1Y,J_3Z)+
\nabla^{CP}F(J_1Y;J_1Z,J_3X)+ \nabla^{CP}F(J_1Z;J_1X,J_3Y).
\end{gather*}
Define $G=F( \cdot,J_2 \cdot)$, we have the following sequence of
equalities
\begin{gather*}
-dG(J_1X,J_1Y,J_1Z)=\\-\nabla^{CP}G(J_1X;J_1Y,J_1Z)-
\nabla^{CP}G(J_1Y;J_1Z,J_1X)- \nabla^{CP}G(J_1Z;J_1X,J_1Y)=\\
=\nabla^{CP}F(J_1X;J_1Y,J_3Z)+ \nabla^{CP}F(J_1Y;J_1Z,J_3X)+
\nabla^{CP}F(J_1Z;J_1X,J_3Y)=\\=-dF(J_3X,J_3Y,J_3Z) \end{gather*}
and therefore $dF(J_3X,J_3Y,J_3Z)=dG(J_1X,J_1Y,J_1Z)$. In other
words,$d_3F_3=d_1F_1$ with $F_3=F$ and $F_1=G$.
\end{proof}
\begin{lem}\label{int3}
A PHKT metric locally admits a potential if and only if the
corresponding HPKT-form is locally $D$-exact.
\end{lem}
\begin{proof}
Suppose that $F=- \frac{1}{2}(dd_3+d_1d_2) \mu$. Then $F= \frac{1}{2}(d
\theta+J_1d \theta)$ with $\theta=-J_3d \mu=-d_3\mu$.
Note that $d\theta$ is $(1,1)$-form (for $J_3$) since $d\theta =-dd_3 \mu$.
Therefore $F=D\theta$ according to \eqref{dhp1}.

Conversely, suppose that $F=D\theta$ for some 1-form $\theta$.
Since F is a $(1,1)$-form for $J_3$, we obtain from \eqref{dhp1}
$$d \theta \in \Lambda^{(1,1)}_{J_3},\quad F= \frac{1}{2}(d
\theta+J_1d \theta).$$ Since $J_3$ is an integrable complex
structure, the local $dd_3$-lemma holds: locally there exists
$\mu$ such that $d\theta =-dd_3 \mu$. We get then $$F=
\frac{1}{2}(d \theta+J_1d \theta)=- \frac{1}{2}(dd_3+d_1d_2)
\mu.$$
\end{proof}
\begin{thm}\label{integ}
On an ($4n\ge 8$)-dimensional manifold any HPKT metric admits
locally an HPKT potential or equivalently any $D$-closed HPKT-form
is locally $D$-exact.
\end{thm}
\begin{proof}
Any ($4n\ge 8$)-dimensional HPKT-manifold is a manifold with a
structure group contained in $Sp(n,\mathbb R)\subset GL(n,\tilde
H)$ and it is $1$-integrable due to the existence of a
torsion-free $GL(n,\tilde H)$-connection, the complex product
connection \cite{An,ASal}. Therefore, since the operator $D$ is
$GL(n,\tilde H)$-invariant, it is sufficient to show that any
$D$-closed HPKT-form is locally $D$-exact in $\mathbb R^{4n}$ with
the standard HPKT structure.

In $\mathbb R^{4n}$ we split the complex coordinates into two sets
$(\{z^j,w^j\}, j=1,\dots,n)$.The hyperparacomplex structure is
given by
\begin{gather*}
J_1=-dz^j\otimes\frac{\partial}{\partial
z^j}-dw^j\otimes\frac{\partial}{\partial w^j}+d\bar
z^j\otimes\frac{\partial}{\partial \bar z^j}+d\bar
w^j\otimes\frac{\partial}{\partial \bar w^j},\\
J_2=id\bar w^j\otimes\frac{\partial}{\partial z^j}-id\bar
z^j\otimes\frac{\partial}{\partial w^j}+idw^j\otimes\frac{\partial}{\partial
\bar z^j}-idz^j\otimes\frac{\partial}{\partial \bar w^j},\\
J_3=id\bar w^j\otimes\frac{\partial}{\partial z^j}-id\bar
z^j\otimes\frac{\partial}{\partial w^j}-idw^j\otimes\frac{\partial}{\partial
\bar z^j}+idz^j\otimes\frac{\partial}{\partial \bar w^j}.
\end{gather*}
Let $F_3\in\Lambda^{1,1}_{J_3}, DF_3=0$. The hyper-parahermitian
condition for the 2-tensor $h=F_3(.,J_3)$ implies that
\begin{equation}\label{c3}
h_{z^j\bar z^k}=h_{w^k\bar w^j}; \quad h_{z^j\bar w^k}=-h_{z^k\bar w^j}.
\end{equation}
The condition $d_1F_1=d_3F_3$ becomes
\begin{gather}\label{c4}
h_{z^j\bar w^k,\bar w^l}+h_{z^k\bar w^l,\bar w^j}+h_{z^l\bar w^j,\bar
w^k}=0,
\quad
h_{w^j\bar z^k,w^l}+h_{w^k\bar z^l,w^j}+h_{w^l\bar z^j,w^k}=0\\ \nonumber
h_{w^j\bar z^k,\bar z^l}+h_{w^k\bar z^l,\bar z^j}+h_{w^l\bar z^j,\bar
z^k}=0,
\quad
h_{z^j\bar w^k,z^l}+h_{z^k\bar w^l,z^j}+h_{z^l\bar w^j,z^k}=0,\\ \nonumber
h_{z^j\bar z^l,\bar w^k}-h_{z^k\bar z^l,\bar w^j}-h_{z^j\bar w^k,\bar
z^l}=0,\quad
h_{z^j\bar z^k,w^l}-h_{z^j\bar z^l,w^k}+h_{w^k\bar z^l,z^j}=0,\\ \nonumber
h_{z^j\bar z^k,\bar z^l}-h_{z^j\bar z^l,\bar z^k}-h_{w^k\bar z^l,\bar
w^j}=0, \quad
h_{z^j\bar z^l,z^k}-h_{z^k\bar z^l,z^j}+h_{z^j\bar w^k,w^l}=0.
\end{gather}
The first and second lines of \eqref{c4}, when combined with the
antisymmetry in $j,k$ of $h_{z^j\bar w^k}$, allow us to apply the
local $\partial\bar{\partial}$-lemma. Therefore, we can write
\begin{equation}\label{c5}
h_{z^j\bar w^k}=(\partial_{z^j}\partial_{\bar w^k}-
\partial_{z^k}\partial_{\bar w^j})\mu; \quad
h_{w^j\bar z^k}=(\partial_{w^j}\partial_{\bar z^k}-
\partial_{w^k}\partial_{\bar z^j})\mu,
\end{equation}
where $\mu$ is some (real) (by para-hermiticity of the metric-and therefore
identical in the two equations \eqref{c5}) function. Inserting \eqref{c5}
into the third equation of \eqref{c4} gives
\begin{equation}\label{c6}
\partial_{\bar w^k}(h_{z^j\bar z^l}-\mu,_{z^j\bar z^l})-
\partial_{\bar w^j}(h_{z^k\bar z^l}-\mu,_{z^k\bar z^l})=0,
\end{equation}
and therefore,
\begin{equation}\label{c7}
h_{z^j\bar z^k}=\mu,_{z^j\bar z^k}+\partial_{\bar w^j}\alpha_{\bar z^k}
\end{equation}
for some integration one form $\alpha_{\bar z^k}$. Combining this
with the fourth equation of \eqref{c4} gives $\alpha_{\bar
z^k}=\mu,_{w^k}$. Thus we get that the function $\mu$ generates
$F_3$. The Lemma~\ref{int3} completes the proof.
\end{proof}
\begin{rem}{\bf A hyperbolic version of Salamon's quaternionic operator \cite{S3}.}
We recall that an \emph{almost paraquaternionic structure} on $M$
is a rank-3 subbundle $P \subset End(TM)$ which is locally spanned
by an almost hyper-paracomplex structure $\tilde H=(J_a)$.
Equivalently, the structure group of $TM$ can be reduced to
$GL(n,\tilde H)Sp(1,\RR)$. A linear connection  on $TM$ is called
\emph{paraquaternionic connection} if it preserves $P$. An almost
paraquaternionic structure is said to be a \emph{paraquaternionic}
if there is a torsion-free paraquaternionic connection. The
paraquaternionic condition controls the Nijenhuis tensors in the
sense that $N(X,Y)(J_a):=N_a$(X,Y) preserves the subbundle $P$.
When $n\ge 2$, the paraquaternionic condition is a strong
condition which is equivalent to the 1-integrability of the
associated $GL(n,\tilde H)Sp(1,\mathbb R)$- structure
\cite{An,ASal}. We can extend the hyper-paracomplex operator $D$
defining it on an almost paraquaternionic manifold locally in the
same way. Consequently, Theorem~\ref{int} is also true, namely an
almost paraquaternionic manifold is paraquaternionic exactly when
$D^2=0$. The proof of Theorem~\ref{int} goes through in this case
also. Now the $Sp(1,\RR)$-part of the paraquaternionic connection
used, adds an additional $Sp(1,\RR)$ term in formula \eqref{new1}
which reflects on \eqref{new2}, whence the Nijenhuis tensors
preserve the subbundle $P$. Using the 1-integrability of the
paraquaternionic structure and the proof of Theorem~\ref{integ}
one gets the local exactness of certain $D$-closed  forms.
\end{rem}

\bibliographystyle{hamsplain}

\providecommand{\bysame}{\leavevmode\hbox to3em{\hrulefill}\thinspace}






\end{document}